\documentclass{article}

\def\nkyear{201?}
\def\nkvol{3?}
\def\nkno{?}

\def\firstpage{1}
\def\lastpage{20}
\def\correction{4}
\def\shorttitle{Homogenization with the quasistatic Tresca friction law} % 压缩标题, 尽量不超过60个字符 (包括空格)
\def\shortauthor{{\it C. Ye,\ \ E. T. Chung\ \ and\ \ J. Cui}} % 作者
\usepackage b \baselineskip 14.5pt

\input{cam.doi}
\journalnumber{11401}
\DOImsnr{0001}
\DOIyear{2007}

\title{\Large \bf \boldmath\ \\ Homogenization with the quasistatic Tresca friction law: qualitative and quantitative results} % 完整标题

\author{\large  Changqing Ye$^1$\qquad Eric T. Chung$^{2}$\qquad Jun-Zhi Cui$^{3}$} % 作者全名

\date{}

\usepackage{mathtools}

\usepackage{mismath}

\usepackage{esint}

\usepackage{subcaption}

\usepackage{siunitx}

\usepackage{xcolor}
\definecolor{GeoDataViz-D1}{HTML}{009392}
\definecolor{GeoDataViz-D2}{HTML}{39B185}
\definecolor{GeoDataViz-D3}{HTML}{9CCB86}
\definecolor{GeoDataViz-D4}{HTML}{E9E29C}
\definecolor{GeoDataViz-D5}{HTML}{EEB479}
\definecolor{GeoDataViz-D6}{HTML}{E88471}
\definecolor{GeoDataViz-D7}{HTML}{CF597E}

\usepackage{tikz}
\usetikzlibrary{calc,patterns}
\tikzstyle{load}   = [thick,-latex]
\tikzstyle{stress} = [-latex]
\tikzstyle{dim}    = [latex-latex]
\tikzstyle{axis}   = [-latex,black!55]
\tikzstyle{label} = [-latex,dashed]

\usepackage{cleveref}
\crefname{lem}{lemma}{lemmas}
\Crefname{lem}{Lemma}{Lemmas}
\crefname{prop}{proposition}{propositions}
\Crefname{prop}{Proposition}{Propositions}
\crefname{thm}{theorem}{theorems}
\Crefname{thm}{Theorem}{Theorems}
\crefname{prob}{problem}{problems}
\Crefname{prob}{Problem}{Problems}
\crefname{defi}{definition}{definitions}
\Crefname{defi}{Definition}{Definitions}
\crefname{equation}{}{}

\newcommand{\Real}{\mathbb{R}}

\newcommand{\Integer}{\mathbb{Z}}

\DeclareMathOperator{\Div}{div}
\DeclareMathOperator{\Curl}{curl}

\DeclarePairedDelimiter{\RoundBrackets}{(}{)}
\DeclarePairedDelimiter{\CurlyBrackets}{\{}{\}}

\DeclarePairedDelimiter{\AngleBrackets}{\langle}{\rangle}

\newcommand{\GammaD}{{\Gamma_{\scriptscriptstyle \mathup{D}}}}
\newcommand{\GammaN}{{\Gamma_{\scriptscriptstyle \mathup{N}}}}
\newcommand{\GammaC}{{\Gamma_{\scriptscriptstyle \mathup{C}}}}
\newcommand{\TrescaBound}{H_{\scriptscriptstyle \mathup{T}}}

\newcommand{\tenA}{\mathsf{A}}

\newcommand{\jrm}{\mathrm{j}}

\newcommand{\Jop}{\mathcal{J}}
\newcommand{\Lop}{\mathcal{L}}
\newcommand{\Mop}{\mathcal{M}}
\newcommand{\Rop}{\mathcal{R}}
\newcommand{\Sop}{\mathcal{S}}

\begin{document}

\maketitle

\thispagestyle{first}
\renewcommand{\thefootnote}{\fnsymbol{footnote}}

% 单位、地址、基金
\footnotetext{\hspace*{-5mm} \begin{tabular}{@{}r@{}p{13.4cm}@{}}
& Manuscript received  \\ %date
$^{1}$ & Department of Mathematics, The Chinese University of Hong Kong,
Shatin, Hong Kong SAR.\\
&{E-mail:cqye@math.cuhk.edu.hk} \\
$^{2}$ & Department of Mathematics, The Chinese University of Hong Kong,
Shatin, Hong Kong SAR.\\
&{E-mail:tschung@math.cuhk.edu.hk} \\
$^{3}$ & LSEC, Academy of Mathematics and Systems Science, Chinese Academy of Sciences,
Beijing 100190, China.\\
&{E-mail:cjz@lsec.cc.ac.cn} \\

$^{\ast}$ & 
The research of Eric Chung is partially supported by the Hong Kong RGC General Research Fund (Project numbers: 14305222 and 14304021). The research of Jun-Zhi Cui is supported by the National Natural Science Foundation of China (No.51739007) and Strategic Priority Research Program of the Chinese Academy of Sciences (No.XDC06030101).
\end{tabular}}

\renewcommand{\thefootnote}{\arabic{footnote}}

\begin{abstract} % 摘要
Modeling of frictional contacts is crucial for investigating mechanical performances of composite materials under varying service environments. The paper considers a linear elasticity system with strongly heterogeneous coefficients and quasistatic Tresca friction law, and studies the homogenization theories under the frameworks of H-convergence and small $\epsilon$-periodicity. The qualitative result is based on H-convergence, which shows the original oscillating solutions will converge weakly to the homogenized solution, while our quantitative result provides an estimate of asymptotic errors in $H^1$-norm for the periodic homogenization. This paper also designs several numerical experiments to validate the convergence rates in the quantitative analysis.

\vskip 4.5mm

\nd \begin{tabular}{@{}l@{ }p{10.1cm}} {\bf Keywords } &
%关键词
Homogenization, Frictional contact mechanics, Quasistatic Tresca friction law
\end{tabular}

\nd {\bf 2000 MR Subject Classification } % 分类号
74M15, 74Q05

\end{abstract}

\baselineskip 14pt

\setlength{\parindent}{1.5em}

\setcounter{section}{0}

\Section{Introduction} \label{sec:Intro}

Modeling physical or engineering problems with Partial Differential Equations (PDE) usually contains three key components: fundamental physical principles, constitutive laws, and Boundary Conditions (BC). The fundamental physical principles are well-established in most scenarios; the relations of different physical quantities are described by the constitutive laws, which appear as coefficients of linear PDEs; constructing appropriate BCs is certainly an application-driven topic, while simplified BCs are commonly adopted to focus on the study of differential operators. This paper concerns a linear elasticity system with \emph{strongly heterogeneous coefficients}---the constitutive law from composite materials and \emph{quasistatic Tresca friction law}---a contact BC.

Composite materials are ubiquitous nowadays, and mathematical treatments for those give birth to homogenization theories (\cite{Milton2002,Torquato2002}). For those composite materials that manifest a significant periodic pattern, periodic homogenization may serve as a proper mathematical tool (\cite{Oleinik1992,Jikov1994,Cioranescu1999,Bensoussan2011}). There are several developments that generalize periodic homogenization---H-convergence (\cite{Tartar2009}), G-convergence (\cite{ChiadoPiat1990,Pankov1997}), two-scale convergence (\cite{Allaire1992,Lukkassen2002}), $\Gamma$-convergence (\cite{DalMaso1993,Braides2002}), Mosco-convergence (\cite{Mosco1994}), Bloch wave spectral analysis (\cite{Conca1997,Allaire1998}), and the periodic unfolding method (\cite{Cioranescu2002,Cioranescu2008}), to name a few, and we adopt the framework of H-convergence in our \emph{qualitative} analysis. A powerful technique in the periodic homogenization is taking formal asymptotic expansions (\cite{Bensoussan2011,Engquist2008}), and estimates of those asymptotic errors are termed as quantitative homogenization theories in the literature (\cite{Avellaneda1987,Gloria2017,Armstrong2019}). A celebrated advancement is the optimal convergence rate estimates in $L^p$- or $H^1$-norm for Dirichlet and Neumann boundary value problems (\cite{Kenig2013,Shen2017,Shen2018}), and some analytical techniques are also applicable in our \emph{quantitative} analysis. For random composite materials, the stochastic homogenization theory can provide a mathematical explanation to average physical properties \cite{Torquato2002}. Recently, Yang and Guan \cite{Yang2022,Yang2022a} proposed a new stochastic homogenization method and related error analysis for the PDEs with random fast oscillation coefficients, which effectively improves the stochastic homogenization theory.

It should be addressed that only obtaining marcoscale or homogenized information is inadequate in some simulations, and we are hence required to dig finescale information ``cheaply''. In the homogenization of small periodicity, we can utilize asymptotic expansions which can be constructed from the homogenized solution and correctors (ref. \cite{Bensoussan2011}) to approximate the original multiscale solution (\cite{Temizer2012,Yang2015}). Beyond periodic homogenization, there are some modern multiscale computational methods, such as generalized multiscale finite element methods (\cite{Efendiev2013,Chung2016,Chung2018}), localized orthogonal decomposition (\cite{Maalqvist2014,Altmann2021,Maalqvist2021}) and generalized finite element methods (\cite{Babuska2011,Babuska2020,Ma2022}), which are able to resolve the finescale information on a coarse mesh. Nevertheless, most of those multiscale computational methods highlight on the heterogeneity in PDE's coefficients, and the numerical error theories of those methods are based on Dirichlet or Neumann boundary value problems, which raises the necessity of challenging these methods via the nonlinearity and nonsmooth from sophisticate BCs.

Before presenting the quasistatic Tresca law, we need to introduce some notations first. Let $\Omega\subset \Real^d$ be a bounded domain with a Lipschitz boundary $\Gamma$, where $d=2$ or $3$. The boundary $\Gamma$ is priorly divided into three mutually disjoint open parts $\GammaD$, $\GammaN$ and $\GammaC$ as $\Gamma=\bar{\GammaD}\cup \bar{\GammaN} \cup \bar{\GammaC}$ to impose different BCs. A contact problem with the quasistatic Tresca friction law (ref. \cite{Han2002,Shillor2004}) reads as
\begin{equation} \label{eq:Tresca strong}
\begin{aligned}
&-\Div \sigma = f, &&\text{in} \  \Omega && (\text{\footnotesize the balance equation}),\\
&\sigma = \tenA:\varepsilon(u), &&\text{in} \  \Omega && (\text{\footnotesize the linear elasticity constitutive law}),\\
&u = 0,  &&\text{on} \  \GammaD && (\text{\footnotesize the homogeneous Dirichlet BC}),\\
&\sigma\nu = t, &&\text{on} \  \GammaN && (\text{\footnotesize the inhomogeneous Neumann BC}),\\
&\left\{
\begin{aligned}
&u_\nu = 0,\ \abs{\sigma_\tau} \leq \TrescaBound\\
&\abs{\sigma_\tau} < \TrescaBound \Rightarrow \dot{u}_\tau=0, \\
&\abs{\sigma_\tau} = \TrescaBound \Rightarrow \exists \lambda>0\ \text{s.t.}\  \sigma_\tau=-\lambda\dot{u}_\tau,
\end{aligned}
\right.
&&\text{on} \  \GammaC && (\text{\footnotesize the quasistatic friction law}),
\end{aligned}
\end{equation}
where $\sigma$ is the stress, $\tenA$ is the fourth-order elasticity tensor, $\varepsilon(u)\coloneqq \big(\nabla u+(\nabla u)^\intercal\big)/2$ is the strain of $u$, $u_\nu$ and $u_\tau$ are defined by the normal-tangential decomposition $u= u_\nu+u_\tau$ on $\GammaC$, $\sigma_\nu$ and $\sigma_\tau$ are defined by $\sigma\nu= \sigma_\nu+\sigma_\tau$ similarly, and $\TrescaBound$ is called the Tresca friction bound (see \cite{Shillor2004}). The strong form \cref{eq:Tresca strong} is equivalent to a variational inequality form (ref. \cite{Kikuchi1988}):
\begin{equation}\label{eq:Tresca weak}
\begin{aligned}
&\int_\Omega \varepsilon(v-\dot{u}(t)):\tenA:\varepsilon(u)\di x + \int_{\GammaC} \TrescaBound \RoundBrackets{\abs{v_\tau}-\abs{\dot{u}_\tau(t)}} \di a \\
&\  \geq \int_\Omega f\cdot \RoundBrackets{v-\dot{u}(t)}\di x + \int_{\GammaN} t\cdot \RoundBrackets{v-\dot{u}(t)} \di a,
\end{aligned}
\end{equation}
for all $v \in V$, where $V$ is the properly defined test function space. The mathematical analysis including solvability and regularity of Tresca friction contact problems has been prepared in the monograph \cite{Han2002}.

In the Tresca friction model, the bound $\TrescaBound$ is independent of $\sigma_\nu$ the normal traction on the contact face and the validity of this assumption has been carefully discussed in \cite{Shillor2004}. In general, the value of $\TrescaBound$ is a rough estimate and could be obtained by physical experiments. Another widely accepted model---Coulomb's law of friction---takes $\sigma_\nu$ into account and reads as $\TrescaBound=\mu_\mathup{f}\abs{\sigma_\nu}$, which means Tresca's law is a simplified one. Although Coulomb's formulation is physically intuitive, it brings a huge challenge to mathematical analysis. One reason is that the bound $\mu_\mathup{f}\abs{\sigma_\nu}$ needs the trace of $\sigma$, which is ill-defined unless $\sigma$ possesses a better regularity than $L^2$. Major breakthroughs for the existence of solutions to Coulomb friction problems have been established in \cite{Eck2005}. Due to the heterogeneity of coefficients, the regularity of solutions to our problems will be much worse, which implies that currently the Tresca friction model is the more attainable one in discussing quantitative homogenization.

There are some studies on the homogenization of nonlinear boundary conditions (see \cite{Yosifian1997,Yosifian1999,Yosifian2001,Liu2008,Ye2021}). To the best of our knowledge, few of them consider \emph{quasistatic} friction contacts, and the quantitative result of asymptotic errors in this paper seems to be novel.

This paper is organized as follows. The general notations and settings will be introduced in \cref{sec:Pre}, which also contains the definition of H-convergence of linear elasticity systems and basic concepts in the periodic homogenization. In \cref{sec:general}, beyond the periodic homogenization, we apply H-convergence to show that original oscillating solutions will converge \emph{weakly} the homogenized solution. \Cref{sec:estimate} aims to quantify asymptotic errors, and we derive an estimate in $H^1$-norm. We design several numerical experiments in \cref{sec:Num} for validating our theoretical results. Some conclusion remarks and further discussions are included in \cref{sec:Conclu}.

\Section{Preliminaries} \label{sec:Pre}

Throughout the article, we denote by $X$ a general real Hilbert space, $X^*$ the dual space of $X$, $\AngleBrackets{\cdot, \cdot}_X$ the duality brackets. If $x_n$ converges weakly to $x$ in a certain Banach space, we will write it as $x_n\rightharpoonup x$ as in \cite{Brezis2011}. For a time interval $I=(0,T)$, vector-valued $p$-integrable and Sobolev spaces such as $L^p(I;X)$ and $W^{1,p}(I;X)$ are defined as usual (see \cite{Gasinski2006}), where $1\leq p \leq \infty$. For brevity, we identify any $v\in W^{1,p}(I;X)$ as an \emph{absolutely continuous} vector-valued function (ref. \cite{Gasinski2006}). The summation convention over repeated indices are summed will be adopted.

We first introduce the following problem:
\begin{prob}[general problem] \label{pro:Abstract}
Find a solution $u:\bar{I}\rightarrow X$, such that for a.e.\ $t\in I$
\[
a\big(u(t), v-\dot{u}(t)\big)+\jrm(v)-\jrm\big(\dot{u}(t)\big) \geq \AngleBrackets{f(t), v-\dot{u}(t)}_X, \ \forall v\in X,
\]
and
\[
u(0) = u_\star.
\]
\end{prob}

The regularity of \cref{pro:Abstract} has been studied in \cite{Han1999,Han2002}, which can be summarized into the following theorem.
\begin{thm}\label{thm:general}
Assume that
\begin{itemize}
\item $a(\cdot, \cdot)$ is a symmetric bilinear form on $X$ with
\[
m \norm{v}_X^2\leq a(v, v), \  \abs{a(w, v)} \leq M \norm{w}_X\norm{v}_X, \  \forall v, w \in X;
\]
\item $\jrm(\cdot)$ is a nonnegative, convex, positively homogeneous\footnote{i.e., $\jrm(av)=\abs{a}\jrm(v)$ for all $a\in \Real$ and $v\in X$.} and Lipschitz continuous functional on $X$;

\item $f \in W^{1,p}(I;X^*)$, where $1\leq p \leq \infty$;

\item $u_\star \in X$ and
\[
a(u_\star,v)+\jrm(v)\geq \AngleBrackets{f(0),v}_X, \  \forall v\in X.
\]
\end{itemize}
Then there exists a unique solution $u\in W^{1,p}(I;X)$ of \cref{pro:Abstract} with estimates
\[
\begin{aligned}
\norm{\dot{u}}_{L^p(I;X)} &\leq C_1\norm{\dot{f}}_{L^p(I;X^*)},\\
\norm{u}_{L^\infty(I;X)} &\leq C_2 \RoundBrackets*{\norm{u_\star}_X+\norm{\dot{f}}_{L^1(I;X^*)}},
\end{aligned}
\]
where $C_1$ and $C_2$ are positive constants depending on $m$ and $M$.
\end{thm}

We assume the elasticity tensor $\tenA$ belongs to the following classes.
\begin{defi}[$\Mop(m, M; \mathcal{O})$ and $\Mop^\mathup{s}(m,M; \mathcal{O})$] \label{def:ten cls}
For given $0<m \leq M< \infty$ and an open set $\mathcal{O}\subset \Real^d$, $\Mop(m, M; \mathcal{O})$ denotes the set of rank-$4$-tensor-valued functions with $\tenA \in \Mop(m, M; \mathcal{O})$ satisfying a.e.\ $x\in \mathcal{O}$,
\[
\begin{aligned}
&\tenA_{ij}^{\alpha \beta}(x)=\tenA_{ji}^{\alpha \beta}(x)=\tenA_{ij}^{\beta\alpha}(x), &&\forall i,j,\alpha,\beta\in\CurlyBrackets{1,\dots,d};\\
&\xi:\tenA(x):\xi \geq m \abs{\xi}^2, &&\forall \xi \in \mathbb{S}^d; \\
&\xi:\tenA^{-1}(x):\xi \geq M^{-1} \abs{\xi}^2, &&\forall \xi \in \mathbb{S}^d;
\end{aligned}
\]
where $\tenA^{-1}(x)$ is the inverse of $\tenA(x)$. Moreover, the function class $\Mop^\mathup{s}(m,M; \mathcal{O})$ denotes for the subset \allowbreak of $\Mop(m,M; \mathcal{O})$ with $\tenA \in \Mop^\mathup{s}(m,M; \mathcal{O})$ satisfying an additional symmetric property---a.e.\ $x\in \mathcal{O}$
\[
\tenA_{ij}^{\alpha \beta}(x)=\tenA^{ij}_{\alpha \beta}(x),\ \forall i,j,\alpha,\beta\in\CurlyBrackets*{1,\dots,d}.
\]

% $\CurlyBrackets*{\tenA_{ij}^{\alpha\beta}}$ are entries of $\tenA$, $\mathbb{S}^d$ is the space of symmetric rank-$2$ tensors (matrices),
\end{defi}

We adopt the algebraic definition of H-convergence from \cite{Gustafsson2007}.
\begin{defi}[H-convergence] \label{def:H conv}
Let $\CurlyBrackets{\tenA^{k'}}\subset \Mop(m,M; \Omega)$ and $\tenA^{\infty} \in \Mop(m',M'; \Omega)$. Then $\tenA^{k'}$ H-converges to $\tenA^{\infty}$ if the following holds true:

Whenever $D^{k'}$, $D^\infty$, $E^{k'}$ and $E^\infty \in L^2(\Omega;\Real^{d\times d})$ satisfy
\begin{itemize}
\item $D^{k'}=\tenA^{k'}:E^{k'}$,
\item $D^{k'}\rightharpoonup D^\infty$ and $E^{k'}\rightharpoonup E^\infty$ in $L^2(\Omega; \Real^{d\times d})$,
\item $\CurlyBrackets{\Div D^{k'}}$ is relatively compact in $H^{-1}(\Omega;\Real^d)$,
\item $\CurlyBrackets{\Curl E^{k'}}$ is relatively compact in $H^{-1}(\Omega; \Real^{d\times d \times d})$;
\end{itemize}
then
\[
D^\infty=\tenA^{\infty}:E^{\infty}.
\]
\end{defi}

Note that we put ``$'$'' over the index $k$ here to distinguish it with spatial indices (e.g. $i,j,\alpha,\beta$) which are always chosen from $1$ to $d$. An immediate result shows that $\Mop^\mathup{s}(m,M; \Omega)$ is sequentially closed w.r.t. H-convergence.
\begin{prop}[see \cite{Gustafsson2007}] \label{prop:H close}
If $\CurlyBrackets{\tenA^{k'}}\subset \Mop^\mathup{s}(m,M; \Omega)$, $\tenA^{\infty} \in \Mop(m',M'; \Omega)$, and $\tenA^{k'}$ H-converges to $\tenA^{\infty}$. Then $\tenA^{\infty} \in \Mop^\mathup{s}(m,M; \Omega)$.
\end{prop}

\begin{defi}[$1$-periodicity]
A scalar/vector/tensor-valued function $f$ is called $1$-periodic if for all $z\in \Integer^d$,
\[
f(x+z)=f(x) \  \text{a.e.} \  x\in \Real^d.
\]
\end{defi}

For smooth $1$-periodic functions, the complements w.r.t. different Sobolev norms are marked with the subscript ``$\#$'', e.g. $H^1_{\#}(Y,\Real^d)$, where $Y$ is always the unit cube $(-1/2, 1/2)^d$ (ref. \cite{Cioranescu1999}).

An essential crux of the periodic homogenization is correctors, which provide an explicit expression of effective coefficient.
\begin{defi}[correctors and effective coefficients from \cite{Shen2018}]
Let $\tenA$ be $1$-periodic and belong to $\Mop^\mathup{s}(m,M;\Real^d)$. The set of correctors of $\tenA$ satisfies $\CurlyBrackets*{\chi_k^\gamma}_{1\leq k,\gamma\leq d}\ \subset H^1_\#(Y;\Real^d)$, where $\chi_k^\gamma$ is the unique solution of the following variational problem:
\[
\left\{
\begin{aligned}
&\chi_k^\gamma = [\chi_k^{1\gamma},\dots,\chi_k^{d\gamma}] \in H^1_\#(Y;\Real^d) \  \text{with} \  \int_Y \chi_k^\gamma(y) \di y=0 \\
&\  \text{and} \  \forall v=[v^1,\dots,v^d] \in H^1_\#(Y;\Real^d), \\
&\int_Y \tenA_{ij}^{\alpha\beta}(y) \partial_j\chi_k^{\beta\gamma}\partial_i v^\alpha \di y = -\int_Y \tenA_{ik}^{\alpha\gamma}(y) \partial_i v^\alpha \di y.
\end{aligned}
\right.
\]
The effective coefficient tensor of $\tenA$ is denoted by $\hat{\tenA}$, and the component $\hat{\tenA}_{ik}^{\alpha\gamma}$ of $\hat{\tenA}$ is defined by
\[
\hat{\tenA}_{ik}^{\alpha\gamma}\coloneqq  \fint_Y \CurlyBrackets*{\tenA_{ik}^{\alpha\gamma}(y)+\tenA_{ij}^{\alpha\beta}\partial_j \chi_k^{\beta\gamma}(y)} \di y.
\]
\end{defi}

H-convergence stems from the periodic homogenization by the following proposition.
\begin{prop}[see \cite{Gustafsson2007}]
Let $\tenA$ belong to $\Mop^\mathup{s}(m,M;\Real^d)$ and also be $1$-periodic. Then $\tenA^\epsilon(x)=\tenA(x/\epsilon)\in \Mop^\mathup{s}(m,M;\Omega)$, and $\tenA^\epsilon$ H-converges to $\hat{\tenA}$ as $\epsilon$ tends to $0$.
\end{prop}

In \cref{sec:general,sec:estimate}, we will depart from concrete quasistatic Tresca friction contact problems as introduced in \cref{sec:Intro} to a more general one. We denote by $V$ a Hilbert space, which is a subspace of $H^1(\Omega;\Real^{d})$ and also satisfies the homogeneous Dirichlet BC and other constraints. For example, $V$ is
\[
\CurlyBrackets*{v \in H^1(\Omega;\Real^d)\mid v = 0 \  \text{a.e.\ on} \  \GammaD,\ \nu\cdot v = v_{\nu} = 0 \  \text{a.e.\ on} \  \GammaC}
\]
for the problem in \cref{sec:Intro}. The notation $V_0$ is equivalent to $H^1_0(\Omega;\Real^d)$, and we \emph{assume} that
\[
V_0 \subset V \subset H^1(\Omega;\Real^d).
\]
Moreover, we \emph{assume} the norm of $V$ (also $H^1(\Omega;\Real^d)$) could be defined by the seminorm of $H^1(\Omega;\Real^d)$ and Korn's inequality always holds (see \cite{Oleinik1992}).

\Section{General results under H-convergence} \label{sec:general}

The model problems we consider in this section are stated as follows.
\begin{prob}\label{pro:H-con k}
Find a solution $u^{k'}:\bar{I}\rightarrow V$, such that for a.e.\ $t\in I$,
\begin{equation} \label{eq:varia k}
a^{k'}\big(u^{k'}(t), v-\dot{u}^{k'}(t)\big)+\Jop(v)-\Jop\big(\dot{u}^{k'}(t)\big) \geq \AngleBrackets*{\Lop(t), v-\dot{u}^{k'}(t)}_V, \  \forall v \in V
\end{equation}
and
\[
u^{k'}(0)=u^{k'}_\star,
\]
where
\[
a^{k'}(w,v)=\int_\Omega \varepsilon(v): \tenA^{k'}:\varepsilon(w) \di x.
\]
\end{prob}

\begin{prob}\label{pro:H-con eff}
Find a solution $u^{\infty}:\bar{I}\rightarrow V$, such that for a.e.\ $t\in I$,
\begin{equation} \label{eq:varia inf}
a^{\infty}\big(u^{\infty}(t), v-\dot{u}^{\infty}(t)\big)+\Jop(v)-\Jop\big(\dot{u}^{\infty}(t)\big) \geq \AngleBrackets*{\Lop(t), v-\dot{u}^{\infty}(t)}_V, \  \forall v \in V
\end{equation}
and
\[
u^{\infty}(0)=u^{\infty}_\star,
\]
%\begin{equation}
%a^{\infty}(u^{\infty}(0), v)+\Jop(v) \geq \AngleBrackets{\Lop(0), v}_V,
%\end{equation}
where
\[
a^{\infty}(w,v)=\int_\Omega \varepsilon(v): \tenA^{\infty}:\varepsilon(w) \di x.
\]
\end{prob}

The following theorem is our main result in this section.
\begin{thm}\label{thm:H main}
Assume that
\begin{itemize}
\item $\CurlyBrackets{\tenA^{k'}}\subset \Mop^\mathup{s}(m,M;\Omega)$ and $\tenA^{k'}$ H-converges to $\tenA^\infty$;
\item $\Jop(\cdot)$ is a nonnegative, convex, positively homogeneous and Lipschitz continuous functional on $V$, and
\begin{equation}\label{eq:ass J}
\Jop(v+w)=\Jop(v), \ \forall v \in V \  \text{and} \  w \in V_0;
\end{equation}
\item $\Lop \in W^{1,p}(I;V^*)$, where $1 < p < \infty$;
\item $u^{k'}_\star$, $u^\infty_\star\in V$ and are respectively the solutions of
%\item $u^{k'}_*\rightharpoonup u^\infty_*$ in $V$, $\Jop(u^{k'}_*)=\Jop(u^\infty_*)=0$, and for all $v\in V$,
\begin{align}
a^{k'}(u^{k'}_\star, v) &= \AngleBrackets{\Lop(0), v}_V, &&\forall v\in V, \label{eq:init k}\\
a^{\infty}(u^{\infty}_\star, v) &= \AngleBrackets{\Lop(0), v}_V, &&\forall v\in V.\label{eq:init inf}
\end{align}
\end{itemize}
Then there exist a unique solution $u^{k'}\in W^{1,p}(I;V)$ of \cref{pro:H-con k} and a unique solution $u^\infty \in W^{1,p}(I;V)$ of \cref{pro:H-con eff}, such that $u^{k'}$ converges weakly to $u^\infty$ in $W^{1,p}(I;V)$.
\end{thm}

\begin{rem}
Comparing with \cref{thm:general}, a notable difference is that we require \cref{eq:ass J}, which essentially says that $\Jop(v)$ is solely determined by the trace of $v$ on $\Gamma$. For the Tresca friction problem stated in \cref{eq:Tresca weak}, we can verify that
\[
\Jop(v) \coloneqq \int_{\GammaC} \TrescaBound \abs{v_\tau} \di a
\]
satisfies all requirements for $\Jop(\cdot)$ in \cref{thm:H main}.
\end{rem}

To prove the main theorem of this section, we need several lemmas prepared. The next lemma shows that the weak convergence on $W^{1,p}(I;V)$ will induce weak convergences on $V$ at any time point.
\begin{lem} \label{lem:e conv}
If $u_n\rightharpoonup u_\infty$ in $W^{1,p}(I;X)$ for $1\leq p \leq \infty$. Then $u_n(t)\rightharpoonup u_\infty(t)$ in $X$ for all $t \in \bar{I}$.
\end{lem}
\begin{proof}
We first prove $u_n(0)\rightharpoonup u_\infty(0)$. For any given $u'\in X^*$, identify it as a functional on $W^{1,p}(I;X)$ via $v \mapsto \AngleBrackets{u',v(0)}_X$, which is a bounded linear map by the embedding $W^{1,p}(I;X)\hookrightarrow C(\bar{I};X)$. Then we have $u_n(0)\rightharpoonup u_\infty(0)$.

For any given $u'\in X^*$ and $t\in \bar{I}$, the linear functional $v\mapsto \AngleBrackets{u',\int_{0}^t\dot{v}(s)\di s}$ is bounded on $W^{1,p}(I;X)$ by
\[
\abs{\AngleBrackets{u',\int_{0}^t\dot{v}(s)\di s}}\leq \norm{u'}_{X^*}\int_{0}^{t}\norm{\dot{v}(s)}_X\di s \leq t^{1-1/p}\norm{u'}_{X^*} \norm{\dot{v}}_{L^p(I;X)}.
\]
Then, applying the fundamental theorem of calculus, we have
\[
\begin{aligned}
\AngleBrackets{u', u_n(t)}_X=&\AngleBrackets{u',u_n(0)}+\AngleBrackets{u', \int_{0}^{t}\dot{u}_n(s)\di s}_X \\
\underset{n}{\longrightarrow}&\AngleBrackets{u',u_\infty(0)}+\AngleBrackets{u', \int_{0}^{t}\dot{u}_\infty(s)\di s}_X= \AngleBrackets{u', u_\infty(t)}_X,
\end{aligned}
\]
which finishes the proof.
\end{proof}

Note that \cref{eq:varia k,eq:varia inf} are stated for a.e.\ $t\in I$, and it is sometimes convenient to consider a ``time-integral'' form. We present the following lemma, which is a supplement of \cref{thm:general}, and its proof could be found in \cite{Han1999,Han2002}.
\begin{lem}\label{lem:integral}
Let the assumptions of \cref{thm:general} be fulfilled and $u$ be the unique solution of \cref{pro:Abstract}. Then for all $w \in L^1(I;X)$,
\[
\int_I a\big(u(s), w(s)-\dot{u}(s)\big)\di s+\int_I \jrm(w(s)) \di s- \int_I \jrm(\dot{u}(s)) \di s \geq \int_I \AngleBrackets{f(s), v(s)-\dot{u}(s)}_X \di s.
\]
\end{lem}

The famous div-curl lemma is required in the proof of \cref{thm:H main}. The original form of the div-curl lemma is presented for vector fields, while we take some modifications for consistency with our elasticity system setting.
\begin{lem}[see \cite{Jikov1994,Tartar2009}] \label{lem:div-curl}
If $D^{k'}$, $D^\infty$, $E^{k'}$ and $E^\infty \in L^2(\Omega;\mathbb{S}^d)$ such that $D^{k'}\rightharpoonup D^\infty$ and $E^{k'}\rightharpoonup E^\infty$ in $L^2(\Omega;\Real^{d\times d})$. Assume that the components of $\Div(D^{k'})$ and $\Curl(E^{k'})$ are all contained in a compact subset of $H^{-1}(\Omega)$. Then for all $\phi\in C^\infty_0(\Omega)$, we have
\[
\int_\Omega D^{k'}:E^{k'} \phi\di x \underset{k'}{\rightarrow} \int_\Omega D^{\infty}:E^{\infty} \phi\di x,
\]
i.e., $D^{k'}:E^{k'}$ converges weakly to $D^{\infty}:E^{\infty}$ as distributions.
\end{lem}

Now we turn to prove \cref{thm:H main}.
\begin{proof}
Because $\Jop$ is a nonnegative functional, the initial conditions $u^{k'}_\star$ and $u^\infty_\star$ satisfy the requirements of \cref{thm:general}. Therefore, the existence and uniqueness of solutions $u^{k'}$ and $u^\infty$ follow from \cref{thm:general}, noting here that the coercivity and symmetry of $\tenA^\infty$ can be shown by \cref{prop:H close}. Moreover, we have a uniform bound for $u^{k'}$ as
\[
\begin{aligned}
\norm{u^{k'}}_{W^{1,p}(I;V)} &\leq C\CurlyBrackets*{\norm{u^{k'}_\star}_V+\norm{\dot{\Lop}}_{L^p(I;V^*)}} \\
&\leq C\CurlyBrackets*{\norm{\Lop(0)}_{V^*}+\norm{\dot{\Lop}}_{L^p(I;V^*)}}<\infty.
\end{aligned}
\]
Because $0<p<1$ and $V$ is a Hilbert space and $W^{1,p}(I;V)$ is a reflexive space. Then it is possible to extract a subsequence of $\CurlyBrackets{k'}$ such that $u^{k'}\rightharpoonup u^o$ in $W^{1,p}(I;V)$ (we still denote the indices of the subsequence by $k'$.), and the proof will be completed by showing that $u^o$ is exactly $u^\infty$.

By \cref{lem:integral}, we have for all $w \in L^1(I;V)$,
\begin{equation} \label{eq:time int varia}
\begin{aligned}
&\int_I a^{k'}\big(u^{k'}(s), w(s)-\dot{u}^{k'}(s)\big) \di s+\int_I \Jop\big(w(s)\big)\di s-\int_I \Jop\big(\dot{u}^{k'}(s)\big)\di s \\
&\qquad \geq \int_I \AngleBrackets{\Lop(s), w(s)-\dot{u}^{k'}(s)}_V\di s.
\end{aligned}
\end{equation}
Because $u^{k'}$ has a regularity of $W^{1,p}$ w.r.t.\ $t$ and the bilinear form $a^{k'}\RoundBrackets{\cdot,\cdot}$ is symmetric also irrelevant with the time variable (ref. \cite{Gasinski2006}), we have
\[
\begin{aligned}
\int_I a^{k'}\big(u^{k'}(s), \dot{u}^{k'}(s)\big) \di s &= \frac{1}{2} \CurlyBrackets*{a^{k'}\big(u^{k'}(T), u^{k'}(T)\big)-a^{k'}\big(u^{k'}(0), u^{k'}(0)\big)} \\
&=\frac{1}{2} \CurlyBrackets*{a^{k'}\big(u^{k'}(T), u^{k'}(T)\big)-a^{k'}\big(u^{k'}_\star, u^{k'}_\star\big)}.
\end{aligned}
\]

We first analyze $a^{k'}\big(u^{k'}(T), u^{k'}(T)\big)$. According to \cref{lem:e conv}, we have $u^{k'}(T)\rightharpoonup u^o(T)$ in $V$. Let $0<\delta<1$ and $v\in V_0$, construct test functions $\dot{u}^{k'}(\cdot)\pm 1_{(T-\delta,T)}(\cdot)v/\delta \in L^1(I;V)$ and take them into \cref{eq:time int varia}. Applying the assumption \cref{eq:ass J}, we can get
\[
a^{k'}\big(\fint_{T-\delta}^T u^{k'}(s)\di s, v\big)=\big<\fint_{T-\delta}^T\Lop(s)\di s, v\big>_V.
\]
Owing to the strong continuity of $u^{k'}$ and $\Lop$ w.r.t. $t$ and taking $\delta\rightarrow 0$, we could recover
\[
a^{k'}\RoundBrackets{u^{k'}(T),v}=\AngleBrackets{\Lop(T), v}_V,
\]
which is equivalent to
\[
-\Div\RoundBrackets*{\tenA^{k'}:\varepsilon(u^{k'}(T))}=\Lop(T)
\]
in the sense of $H^{-1}$. Denoted by $E^{k'}\coloneqq \varepsilon(u^{k'}(T))$ and $D^{k'}\coloneqq \tenA^{k'}:E^{k'}$, we can see that $E^{k'} \rightharpoonup E^o\coloneqq \varepsilon(u^o(T))$ in $L^2(\Omega;\Real^{d\times d})$ due to $u^{k'}(T) \rightharpoonup u^o(T)$ in $V$. Note that $\CurlyBrackets{D^{k'}}$ is a bounded subset of $L^2(\Omega;\Real^{d\times d})$, therefore we can remove a subsequence of $\CurlyBrackets{k'}$ such that $D^{k'} \rightharpoonup D^0$. Applying \cref{def:H conv} for $\tenA^{k'}$, we deduce that $D^o=\tenA^\infty:E^o$. Moreover, combining \cref{lem:div-curl}, we can also see that $D^{k'}:E^{k'}$ converges weakly to $D^o:E^o$ as distributions. Considering $D^{k'}:E^{k'}=E^{k'}:\tenA:E^{k'}$ is nonnegative a.e.\ in $\Omega$, arbitrarily choosing $\phi \in C^\infty_0(\Omega)$ and $0\leq \phi(x) \leq 1$ for all $x\in \Omega$, we can show the following relation:
\[
a^{k'}\big(u^{k'}(T),u^{k'}(T)\big) \geq \int_\Omega D^{k'}:E^{k'} \phi \di x \underset{k'}{\rightarrow} \int_\Omega D^o:E^o \phi \di x,
\]
which gives that
\[
\liminf_{k'} a^{k'}\big(u^{k'}(T),u^{k'}(T)\big) \geq a^\infty\RoundBrackets{u^o(T),u^o(T)},
\]
up to a subsequence of $\CurlyBrackets{k'}$.

We then handle $a^{k'}(u^{k'}_\star, u_\star^{k'})$. Similarly, we have $u^{k'}(0)=u^{k'}_\star\rightharpoonup u^o(0)$ in $V$ by \cref{lem:e conv}. Utilizing the definition of H-convergence again, we can obtain that $\tenA^{k'}:\varepsilon(u^{k'}_\star) \rightharpoonup \tenA^\infty:\varepsilon(u^o(0))$ in $L^2(\Omega;\Real^{d\times d})$ from the variational form \cref{eq:init k}. Take the LHS limit of \cref{eq:init k}, and it shows that for all $v\in V$,
\[
a^\infty\RoundBrackets{u^o(0), v} = \AngleBrackets{\Lop(0), v}_V,
\]
which implies $u^o(0)=u^\infty_\star$ since $u^\infty_\star$ satisfies the same variational form \cref{eq:init inf}. Moreover, replacing $v$ with $u^{k'}$ in \cref{eq:init k} and taking a limit, we arrive at
%
%Take $v=\pm u^{k'}_*$ in \cref{eq:init k} and recall $\Jop(u^{k'}_*)=\Jop(-u^{k'}_*)=0$, and it is obvious to show that
%\[
%a^{k'}(u^{k'}_*,u^{k'}_*) = \AngleBrackets{\Lop(0), u^{k'}_*}.
%\]
%Similarly, we have
%\[
%a^{\infty}(u^{\infty}_*,u^{\infty}_*) = \AngleBrackets{\Lop(0), u^{\infty}_*}.
%\]
%from \cref{eq:init inf}. Recall that $u^{k'}\rightharpoonup u^o$ in $W^{1,p}(I;V)$, this leads $u^{k'}(0)\rightharpoonup u^o(0)$ in $V$ by \cref{lem:e conv}. Combining with the weak convergence assumption for $u^{k'}_*$, we arrive at $u^o(0)=u^\infty_*$ and
\[
a^{k'}(u^{k'}_\star,u^{k'}_\star)=\AngleBrackets{\Lop(0),u^{k'}_\star} \underset{k'}{\rightarrow} \AngleBrackets{\Lop(0), u^\infty_\star}=a^\infty(u^\infty_\star, u^\infty_\star)=a^\infty(u^o(0), u^o(0)).
\]

Through a similar procedure for $s=T$, we can show that for all $s\in I$, $\tenA^{k'}:\varepsilon(u^{k'}(s))\rightharpoonup \tenA^\infty:\varepsilon(u^o(s))$ in $L^2(\Omega;\Real^{d\times d})$. Then combining the embedding relation $W^{1,p}(I;V)\hookrightarrow L^\infty(I;V)$ and the uniform coercivity of $\tenA^{k'}$ (positive constants $m$ and $M$ of \cref{def:ten cls}), we conclude that
\[
\int_I a^{k'}\RoundBrackets{u^{k'}(s),w(s)}\di s \underset{k'}{\rightarrow} \int_I a^\infty\RoundBrackets{u^o(s),w(s)} \di s
\]
by the dominance convergence theorem.

By the assumption for $\Jop$, it is clear to see that $u\mapsto \int_I \Jop\big(\dot{u}(s)\big)\di s$ is a nonnegative, \emph{convex}, positively homogeneous and \emph{Lipschitz continuous} functional on $W^{1,p}(I;
V)$. We hence have
\[
\liminf_{k'} \int_I \Jop\big(\dot{u}^{k'}(s)\big)\di s \geq \int_I \Jop\big(\dot{u}^o(s)\big) \di s.
\]

Note that
\[
\int_I \AngleBrackets{\Lop(s), \dot{u}^{k'}(s)}_V\di s \underset{k'}{\rightarrow} \int_I \AngleBrackets{\Lop(s), \dot{u}^o(s)}_V\di s
\]
follows from the weak convergence of $u^{k'}$ in $W^{1,p}(I;V)$.

Finally, taking the limit superior of the LHS and the limit of the RHS in \cref{eq:time int varia}, up to a subsequence of $\CurlyBrackets{k'}$, we can show that
\[
\begin{aligned}
&\frac{1}{2}\CurlyBrackets*{a^\infty\big(u^o(0),u^o(0)\big)-a^\infty\big(u^o(T),u^o(T)\big)}+\int_I a^\infty\big(u^o(s), w(s)\big)\di s\\
&\  + \int_I \Jop\big(w(s)\big) \di s-\int_I \Jop\big(\dot{u}^o(s)\big) \di s \\
=&\int_I a^\infty\big(u^o(s), w(s)-\dot{u}^o(s)\big)\di s+ \int_I \Jop\big(w(s)\big) \di s-\int_I \Jop\big(\dot{u}^o(s)\big) \di s \\
\geq & \int_I \AngleBrackets{\Lop(s), w(s)-\dot{u}^o(s)}_V\di s, \  \forall w \in L^1(I;V).
\end{aligned}
\]
We can also rewrite this ``time-integral'' variational inequality into an ``a.e.\ $t$'' form \cite{Han2002} and recall that $u^o(0)=u^\infty(0)=u^\infty_\star$, which exactly says that $u^o$ is again a solution of \cref{pro:H-con eff}, then we prove $u^o=u^\infty$.
\end{proof}

\Section{Estimates for the homogenization of small periodicity} \label{sec:estimate}

For brevity, we will abbreviate $H^1(\Omega;\Real^d)$ as $H^1$ and $H^2(\Omega;\Real^d)$ as $H^2$ in this section. The model problems of this section are as follows.

\begin{prob}\label{pro:per}
Find a solution $u_\epsilon:\bar{I}\rightarrow V$, such that for a.e.\ $t\in I$,
\begin{equation} \label{eq:varia per}
a_\epsilon\big(u_\epsilon(t), v-\dot{u}_\epsilon(t)\big)+\Jop(v)-\Jop\big(\dot{u}_\epsilon(t)\big) \geq \AngleBrackets{\Lop(t), v-\dot{u}_\epsilon(t)}_V, \  \forall v \in V
\end{equation}
and
\[
u_\epsilon(0)=u_{\epsilon,\star},
\]
where
\[
a_\epsilon(w,v)=\int_\Omega \varepsilon(v): \tenA\RoundBrackets*{x/\epsilon}:\varepsilon(w) \di x.
\]
\end{prob}

\begin{prob}\label{pro:homo}
Find a solution $u_0:\bar{I}\rightarrow V$, such that for a.e.\ $t\in I$,
\begin{equation} \label{eq:varia homo}
\hat{a}\big(u_0(t), v-\dot{u}_0(t)\big)+\Jop(v)-\Jop\big(\dot{u}_0(t)\big) \geq \AngleBrackets{\Lop(t), v-\dot{u}_0(t)}_V, \  \forall v \in V
\end{equation}
and
\[
u_0(0)=u_{0,\star},
\]
%\begin{equation}
%a^{\infty}(u^{\infty}(0), v)+\Jop(v) \geq \AngleBrackets{\Lop(0), v}_V,
%\end{equation}
where
\[
\hat{a}(w,v)=\int_\Omega \varepsilon(v): \hat{\tenA}:\varepsilon(w) \di x.
\]
\end{prob}

The next theorem is the main result of this section.
\begin{thm} \label{thm:main esti}
Assume that
\begin{itemize}
\item $\tenA \in \Mop^\textnormal{s}(m,M;\Omega)$ and is $1$-periodic, the correctors of $\tenA$ are $\CurlyBrackets{\chi_k^\gamma}$, and the effective coefficient tensor of $\tenA$ is $\hat{\tenA}$;
\item $\Jop(\cdot)$ is a nonnegative, convex, positively homogeneous and Lipschitz continuous functional on $V$, and
\begin{equation}\label{eq:ass J per}
\Jop(v+w)=\Jop(v), \  \forall v \in V \  \text{and} \  w \in V_0;
\end{equation}
\item $\Lop \in H^1(I;V^*)$;
\item $u_{\epsilon,\star}, u_{0,\star}\in V$ and are respectively the solutions of
\begin{align}
a_\epsilon(u_{\epsilon,\star}, v) &= \AngleBrackets{\Lop(0), v}_V, &&\forall v\in V, \label{eq:init k per}\\
\hat{a}(u_{0,\star}, v) &= \AngleBrackets{\Lop(0), v}_V, &&\forall v\in V.\label{eq:init inf per}
\end{align}
\end{itemize}
Let $u_\epsilon$ be the unique solution of \cref{pro:per}, $u_{0}$ be the unique solution of \cref{pro:homo}, denote the residual by
\begin{equation}
r_\epsilon^\alpha=u_\epsilon^\alpha-u_0^\alpha-\epsilon\chi_j^{\alpha\beta}\RoundBrackets*{x/\epsilon} \partial_j u_0^\beta, \  r_\epsilon=[r_\epsilon^1,\dots,r_\epsilon^d],
\end{equation}
provided that there ere the regularities $u_{0}\in H^1\RoundBrackets*{I;H^2}$ and $\CurlyBrackets{\chi_k^\gamma} \subset W^{1,\infty}_\#(Y;\Real^d)$, then
\[
\max_{t\in \bar{I}} \norm{r_\epsilon(t)}_{H^1} \leq C \epsilon^{1/2}\norm{u_{0}}_{H^1\RoundBrackets*{I;H^2}},
\]
where the positive constant $C$ only depends on $m$, $M$, $T$, $\Omega$ and $\max_{k,\gamma}\norm{\chi_k^\gamma}_{W^{1,\infty}_\#(Y;\Real^d)}$.
\end{thm}
\begin{rem}
The expression $u_{0}^\alpha+\epsilon\chi_j^{\alpha\beta}(x/\epsilon)\partial_ju_0^\beta$ is usually called the first-order expansion of $u_\epsilon^\alpha$ in the literature, and certainly we could expand $u_\epsilon$ to higher orders (see \cite{Bensoussan2011}). However, the mismatch on the boundary put a ceiling of convergence rates of asymptotic expansions w.r.t. $\epsilon$. In our problem, low regularities from variational inequality (that is, the homogenized solution $u_0$ may not be smooth even providing sufficiently smooth data) limit the usage of higher order expansions, which require taking high order derivatives of $u_{0}$.

Note that there are no regularity assumptions for the homogenized solution $u^\infty$ in \cref{thm:H main} and $u^\infty \in W^{1,p}(I;V)$ is a direct result of \cref{thm:general}. Nevertheless, we require a regularity of $H^1(I;H^2)$ for $u_0$ to derive a convergence rate estimate, and we skip technical discussions of what conditions that can induce $u_0 \in H^1(I;H^2)$. Referring to \cref{eq:Tresca weak}, we may see that $\TrescaBound$ could be a function of $\GammaC$, and therefore the smoothness of $\TrescaBound$ could affect the regularity of $u_0$. We emphasize that studying the higher regularity than $H^1$ for the solution to a variational inequality is complicated and beyond the scope of the article.
\end{rem}

The next lemma is the main tool in our analysis.
\begin{lem} \label{lem:boud ly esti}
Let $\varphi_\epsilon \in H^1(\Omega;\Real^d)$ and $\varphi_0 \in H^2(\Omega;\Real^d)$, denote by
\[
w_\epsilon^\alpha = \varphi^\alpha_\epsilon-\varphi^\alpha_0-\epsilon \chi_j^{\alpha\beta}\RoundBrackets*{x/\epsilon} \theta_\epsilon \Sop_\epsilon^2 \big(\partial_j \varphi_0^\beta\big), \  w_\epsilon=[w_\epsilon^{1}, \dots, w_\epsilon^{d}],
\]
where $\theta_\epsilon$ is a cutoff function and $\Sop_\epsilon$ is the $\epsilon$ smoothing operator (see \cite{Shen2018}). Then for all $v \in H^1$, we have
\begin{equation}\label{eq:key estimate}
a_\epsilon\RoundBrackets{w_\epsilon, v}=a_\epsilon\RoundBrackets{\varphi_\epsilon,v}-\hat{a}\RoundBrackets{\varphi_0,v}+\AngleBrackets{\Rop_\epsilon \varphi_0,v}_{H^1},
\end{equation}
where $\Rop_\epsilon: H^2\rightarrow (H^1)^*$ is a linear map. Moreover, $\Rop_\epsilon$ is bounded and satisfies
\[
\norm{\Rop_\epsilon}_{\textnormal{op}} \coloneqq \sup_{\substack{\varphi \in H^2\\ \varphi\neq 0}} \sup_{\substack{v\in H^1\\ v\neq 0}} \frac{\abs{\AngleBrackets{\Rop_\epsilon \varphi, v}_{H^1}}}{\norm{\varphi}_{H^2} \norm{v}_{H^1}} \leq C\epsilon^{1/2},
\]
where the positive constant $C$ only depends on $m$, $M$ and $\Omega$.
\end{lem}

Comparing with the first-order expansion in \cref{thm:main esti},  we introduce the cut-off function $\theta_\epsilon$ and $\epsilon$-smoothing operator $\Sop_\epsilon$ here. The reason is that if we assume only $\chi_j^{\alpha\beta}\in H^1_\#(Y)$ and $\phi_0^\beta \in H^2(\Omega)$, $\epsilon \chi_j^{\alpha\beta}(x/\epsilon) \partial_j\phi_0^\beta$ generally does not have a regularity of $H^1$, while $\epsilon \chi_j^{\alpha\beta}(x/\epsilon) \theta_\epsilon \Sop_\epsilon^2\RoundBrackets{\partial_j\phi_0^\beta}$ will always even belong to $H^1_0$. This property will be repeatedly utilized in the following proofs.

\begin{proof}[Proof of \cref{lem:boud ly esti}]
Taking a direct calculation, we have
\[
\begin{aligned}
\partial_m w_\epsilon^\alpha =& \partial_m \varphi_\epsilon^\alpha-\partial_m \varphi_0^\alpha -\RoundBrackets*{\partial_m \chi_j^{\alpha\beta}}(x/\epsilon)\theta_\epsilon \Sop_\epsilon^2\big(\partial_j\varphi_0^\beta\big) \\
&+\epsilon \chi_j^{\alpha\beta}(x/\epsilon) \partial_m \theta_\epsilon \Sop_\epsilon^2\big(\partial_j\varphi_0^\beta\big) -\epsilon \chi_j^{\alpha\beta}(x/\epsilon) \theta_\epsilon \Sop_\epsilon^2\big(\partial_{jm}\varphi_0^\beta\big) \\
=& \partial_m \varphi_\epsilon^\alpha - \theta_\epsilon \Sop_\epsilon^2\big(\partial_m\varphi_0^\alpha\big) - \RoundBrackets*{\partial_m \chi_j^{\alpha\beta}}(x/\epsilon)\theta_\epsilon \Sop_\epsilon^2\big(\partial_j\varphi_0^\beta\big) \\
&+ \theta_\epsilon \Sop_\epsilon^2\big(\partial_m\varphi_0^\alpha\big)-\partial_m\varphi_0^\alpha \\
&+\epsilon \chi_j^{\alpha\beta}(x/\epsilon) \partial_m \theta_\epsilon \Sop_\epsilon^2\big(\partial_j\varphi_0^\beta\big) -\epsilon \chi_j^{\alpha\beta}(x/\epsilon) \theta_\epsilon \Sop_\epsilon^2\big(\partial_{jm}\varphi_0^\beta\big) \\
\coloneqq & \partial_m \varphi_\epsilon^\alpha - \theta_\epsilon \Sop_\epsilon^2\big(\partial_m\varphi_0^\alpha\big) - \RoundBrackets*{\partial_m \chi_j^{\alpha\beta}}(x/\epsilon)\theta_\epsilon \Sop_\epsilon^2\big(\partial_j\varphi_0^\beta\big) + J_m^\alpha.
\end{aligned}
\]
Then
\[
\begin{aligned}
a_\epsilon(w_\epsilon, v)=&\int_\Omega \tenA_{mn}^{\alpha \gamma} (x/\epsilon) \partial_m w_\epsilon^\alpha \partial_n v^\gamma \di x\\
=& \int_\Omega \tenA_{mn}^{\alpha \gamma}(x/\epsilon)\partial_m \varphi_\epsilon^\alpha\partial_n v^\gamma \di x - \int_\Omega \hat{\tenA}_{mn}^{\alpha\gamma} \partial_m \varphi_0^\alpha\partial_n v^\gamma \di x \\
&+ \int_\Omega \hat{\tenA}_{mn}^{\alpha\gamma} \partial_m \varphi_0^\alpha\partial_n v^\gamma \di x - \int_\Omega \hat{\tenA}_{mn}^{\alpha\gamma} \theta_\epsilon \Sop_\epsilon^2\big(\partial_m \varphi_0^\alpha\big)\partial_n v^\gamma \di x \\
&+ \int_\Omega \CurlyBrackets*{\hat{\tenA}_{jn}^{\beta\gamma}-\RoundBrackets*{\tenA_{jn}^{\beta\gamma}+\tenA_{mn}^{\alpha\gamma}\partial_m \chi_j^{\alpha\beta}}(x/\epsilon)}\Sop_\epsilon^2\big(\partial_j\varphi_0^\beta\big) \partial_n v^\gamma \di x \\
&+ \int_\Omega \tenA_{mn}^{\alpha \gamma}  J_m^\alpha \partial_n v^\gamma \di x \\
\coloneqq & a_\epsilon(\varphi_\epsilon,v)-\hat{a}(\varphi_0,v) + K_1+K_2+K_3.
\end{aligned}
\]
Using the estimates in Chap. 3.1 of \cite{Shen2018}, we could derive an estimate of $K_1+K_3$ as
\[
\abs{K_1}+\abs{K_3} \leq C\epsilon^{1/2} \norm{\varphi_0}_{H^2} \abs{v}_{H^1}.
\]
The estimate of $\abs{K_2}\leq C\epsilon^{1/2} \norm{\varphi_0}_{H^2} \abs{v}_{H^1} $ needs the technique of flux correctors, which is also introduced in Chap. 3.1 of \cite{Shen2018}.
\end{proof}

The next lemma presents the convergence rate of initial conditions.
\begin{lem} \label{lem:init esti}
Under the same assumptions of \cref{thm:main esti}, let $w_\epsilon$ be defined by
\[
w^\alpha_\epsilon = u_{\epsilon,\star}^\alpha-u_{0,\star}^\alpha-\epsilon \chi_j^{\alpha\beta}\RoundBrackets*{x/\epsilon} \theta_\epsilon \Sop_\epsilon^2 \big(\partial_j u_{0,\star}^\beta\big)
\]
as in \cref{lem:boud ly esti}. If $u_{0,\star}\in H^2$, then
\[
\norm{w_\epsilon}_{H^1} \leq C\epsilon^{1/2}\norm{u_{0,\star}}_{H^2},
\]
where the positive constant $C$ only depends on $m$, $M$ and $\Omega$.
\end{lem}
\begin{proof}
Because $\theta_\epsilon \in C_0^\infty(\Omega)$ and $\Sop_\epsilon^2\big(\partial_j u_{0,\star}^\beta\big)$ is smooth on $\Omega$, we can see that
\[
\epsilon \chi_j^{\alpha\beta}\RoundBrackets{x/\epsilon} \theta_\epsilon \Sop_\epsilon^2\big(\partial_j u_{0,\star}^\beta\big) \in H^1_0(\Omega),
\]
which leads $w_\epsilon \in V$ due to $V_0\subset V\subset H^1$. Then, replacing $v$ with $w_\epsilon$ in \cref{eq:init k per,eq:init inf per} respectively, we obtain
\[
a_\epsilon\RoundBrackets{u_{\epsilon,\star},w_\epsilon}=\AngleBrackets{\Lop(0),w_\epsilon}_V=\hat{a}\RoundBrackets{u_{0,\star},w_\epsilon}.
\]
Utilizing the estimate \cref{eq:key estimate}, we can show
\[
\begin{aligned}
c\norm{w_\epsilon}_{H^1} \leq a_\epsilon\RoundBrackets{w_\epsilon,w_\epsilon} =& a_\epsilon\RoundBrackets{u_{\epsilon,\star},w_\epsilon}-\hat{a}\RoundBrackets{u_{0,\star},w_\epsilon}+\AngleBrackets{\Rop_\epsilon u_{0,\star},w_\epsilon}_{H^1} \\
=& \AngleBrackets{\Lop(0),w_\epsilon}_V-\AngleBrackets{\Lop(0),w_\epsilon}_V+\AngleBrackets{\Rop_\epsilon u_{0,\star},w_\epsilon}_{H^1} \\
\leq & \norm{\Rop_\epsilon}_{\textnormal{op}} \norm{u_{0,\star}}_{H^2} \norm{w_\epsilon}_{H^1} \\
\leq & C\epsilon^{1/2} \norm{u_{0,\star}}_{H^2} \norm{w_\epsilon}_{H^1},
\end{aligned}
\]
where the positive constant $c$ is from Korn's inequality, and consequently it gives us the target estimate.
\end{proof}

The next lemma paves a path to the proof of \cref{thm:main esti}.
\begin{lem} \label{lem:mid esti}
Under the same assumptions and notations of \cref{thm:main esti}, let $w_\epsilon$ be defined by

\[
w_\epsilon^\alpha = u_{\epsilon}^\alpha-u_0-\epsilon \chi_j^{\alpha\beta}\RoundBrackets*{x/\epsilon} \theta_\epsilon \Sop_\epsilon^2\big(\partial_j u_{0}^\beta\big)
\]
as in \cref{lem:boud ly esti}. If $u_{0}\in H^1(I;H^2)$, then
\[
\max_{t\in \bar{I}} \norm{w_\epsilon(t)}_{H^1} \leq C\epsilon^{1/2} \norm{u_{0}}_{H^1(I;H^2)},
\]
where the positive constant $C$ only depends on $m$, $M$, $T$ and $\Omega$.
\end{lem}

To prove this lemma, we need the following proposition which is essentially Integration by Parts for the vector-valued Sobolev spaces.
\begin{prop}[see Proposition 2.2.43 of \cite{Gasinski2006}] \label{lem:int by part}
Let $u \in H^1(I;X^*)$ and $v \in H^1(I;X)$, we have
\[
\AngleBrackets{u(t),v(t)}_X-\AngleBrackets{u(s),v(s)}_X = \int_{s}^{t} \AngleBrackets{\dot{u}(t), v(t)}_X+\AngleBrackets{u(t),\dot{v}(t)}_X \di t
\]
for all $0\leq s \leq t \leq T$.
\end{prop}
\begin{proof}[Proof of \cref{lem:mid esti}]
We define $u_{1,\epsilon}$ by
\[
u_{1,\epsilon}^\alpha \coloneqq \epsilon \chi_j^{\alpha\beta}\RoundBrackets*{x/\epsilon} \theta_\epsilon \Sop_\epsilon^2\big(\partial_j u_{0}^\beta\big),
\]
which leads $w_\epsilon=u_\epsilon-u_{0}-u_{1,\epsilon}$. Furthermore, we have $u_{1,\epsilon}(t) \in V_0$ for all $t\in \bar{I}$ and $\dot{u}_{1,\epsilon}(t) \in V_0$ for a.e.\ $t\in I$ via the properties of $\theta_\epsilon$ and $\Sop_\epsilon$. Note that $\dot{u}_0(t)+\dot{u}_{1,\epsilon}(t) \in V$, we can hence replace $v$ with it in \cref{eq:varia per}, which implies for a.e.\ $t\in I$,
\[
a_\epsilon\big(u_\epsilon(t), -\dot{w}_\epsilon(t)\big) + \Jop\big(\dot{u}_0(t)+\dot{u}_{1,\epsilon}(t)\big)-\Jop\big(\dot{u}_\epsilon(t)\big) \geq \AngleBrackets{\Lop(t), -\dot{w}_\epsilon(t)}_V.
\]
Note the assumption for $\Jop$ in \cref{eq:ass J per}, then $\Jop\big(\dot{u}_0(t)+\dot{u}_{1,\epsilon}(t)\big)=\Jop\big(\dot{u}_0(t)\big)$ and for a.e.\ $t\in I$
\begin{equation} \label{eq:temp per 1}
a_\epsilon\big(u_\epsilon(t), \dot{w}_\epsilon(t)) + \Jop\big(\dot{u}_\epsilon(t)\big)-\Jop\big(\dot{u}_0(t)\big) \leq \AngleBrackets{\Lop(t), \dot{w}_\epsilon(t)}_V.
\end{equation}
Similarly, it also holds that $\dot{u}_\epsilon(t)-\dot{u}_{1,\epsilon}(t) \in V$ for a.e.\ $t\in V$, and taking it into the variational form \cref{eq:varia homo} gives us that for a.e.\ $t\in I$,
\begin{equation}\label{eq:temp per 2}
\hat{a}\big(u_{0}(t), \dot{w}_\epsilon(t)\big)+\Jop\big(\dot{u}_\epsilon(t)\big)-\Jop\big(\dot{u}_0(t)\big) \geq \AngleBrackets{\Lop(t), \dot{w}_\epsilon(t)},
\end{equation}
where the relation $\Jop\big(\dot{u}_\epsilon(t)-\dot{u}_{1,\epsilon}(t)\big)=\Jop\big(\dot{u}_\epsilon(t)\big)$ is considered. Combining \cref{eq:temp per 1,eq:temp per 2}, we derive an important inequality---
\[
a_\epsilon\big(u_\epsilon(t),\dot{w}_\epsilon(t)\big)-\hat{a}\big(u_{0}(t),\dot{w}_\epsilon(t)) \leq 0
\]
for a.e.\ $t\in I$.

Observing \cref{eq:key estimate} and substituting $\dot{w}(s)$ for $v$, we obtain that for a.e.\ $s\in I$,
\begin{equation}\label{eq:temp per 3}
\begin{aligned}
a_\epsilon\big(w_\epsilon(s), \dot{w}_\epsilon(s)\big) &= a_\epsilon\big( u_\epsilon(s), \dot{w}_\epsilon(s)\big)-\hat{a}\big(u_{0}(s), \dot{w}_\epsilon(s)\big)+\AngleBrackets{\Rop_\epsilon u_{0}(s),\dot{w}_\epsilon(s)}_{H^1} \\
&\leq \AngleBrackets{\Rop_\epsilon u_{0}(s),\dot{w}_\epsilon(s)}_{H^1}.
\end{aligned}
\end{equation}
By the assumption that $\Lop \in H^1(I;V^*)$, we can get $u_\epsilon \in H^1(I;V)$. Moreover, from the regularity $u_{0} \in H^1(I;H^2)$, it is not hard to deduce that $u_{1,\epsilon}\in H^1(I;H^1)$, which shows $w_\epsilon \in H^1(I;H^1)$. According to the linearity of the operator $\Rop_\epsilon$, we naturally obtain that for a.e.\ $s \in I$,
\[
\frac{\di \Rop_\epsilon u_{0}}{\di t}(s)= \Rop_\epsilon \dot{u}_0(s)
\]
in $(H^1)^*$ and $\Rop_\epsilon u_0 \in H^1(I;(H^1)^*)$. Then $\Rop_\epsilon u_{0}$ and $w_\epsilon$ fulfill the requirements of \cref{lem:int by part}, and we hence have
\[
\begin{aligned}
&\int_{0}^t \AngleBrackets{\Rop_\epsilon u_{0}(s),\dot{w}_\epsilon(s)}_{H^1}\di s \\
=&\AngleBrackets{\Rop_\epsilon u_0(t), w_\epsilon(t)}_{H^1}- \AngleBrackets{\Rop_\epsilon u_{0,\star}, w_\epsilon(0)}_{H^1}-\int_0^t \AngleBrackets{\Rop_\epsilon \dot{u}_{0}(s),w_\epsilon(s)}_{H^1}\di s,
\end{aligned}
\]
which gives an expression of the time integration of the RHS of \cref{eq:temp per 3}. Meanwhile, take the time integration of the LHS of \cref{eq:temp per 3}, we have
\[
\int_{0}^t a_\epsilon\big(w_\epsilon(s), \dot{w}_\epsilon(s)\big) \di s= \frac{1}{2} \CurlyBrackets*{a_\epsilon\big(w_\epsilon(t), w_\epsilon(t)\big)-a_\epsilon\big(w_\epsilon(0), w_\epsilon(0)\big)}.
\]
Combining the coercivity of $a_\epsilon\RoundBrackets{\cdot,\cdot}$ and H\"{o}lder's inequality, we obtain
\begin{equation}\label{eq:temp per 4}
\begin{aligned}
\norm{w_\epsilon(t)}_{H^1}^2 \leq & C \Big\{\norm{w_\epsilon(0)}_{H^1}^2+\norm{\Rop_\epsilon u_0(t)}_{(H^1)^*}^2+\norm{\Rop_\epsilon u_{0,\star}}_{(H^1)^*}\norm{w_\epsilon(0)}_{H^1}\\
&\quad +\int_0^t \norm{\Rop_\epsilon \dot{u}_0(s)}_{(H^1)^*}^2 +\norm{w_\epsilon(s)}_{H^1}^2\di s\Big\},
\end{aligned}
\end{equation}
where the positive constant $C$ only depends on $m$, $M$ and $d$.

It has been proved that $\norm{w_\epsilon(0)}_{H^1} \leq C\epsilon^{1/2} \norm{u_{0,\star}}_{H^2}$ in \cref{lem:init esti}. Together with the estimate $\norm{\Rop}_{\textnormal{op}}\leq C\epsilon^{1/2}$ from \cref{lem:boud ly esti}, we can obtain the following estimates:
\[
\begin{aligned}
\norm{w_\epsilon(0)}_{H^1}^2 &\leq C\epsilon \norm{u_{0,\star}}_{H^2}^2,\\
\norm{\Rop_\epsilon u_0(t)}_{(H^1)^*}^2 &\leq \norm{\Rop_\epsilon}_{\textnormal{op}}^2 \norm{u_{0}(t)}^2_{H^2} \leq C\epsilon \norm{u_{0}}_{L^\infty(I;H^2)}^2, \\
\norm{\Rop_\epsilon u_{0,\star}}_{(H^1)^*}\norm{w_\epsilon(0)}_{H^1} &\leq C\epsilon^{1/2}\norm{\Rop_\epsilon}_{\textnormal{op}} \norm{u_{0,\star}}_{H^2}^2 \leq C\epsilon\norm{u_{0,\star}}_{H^2}^2, \\
\int_0^t \norm{\Rop_\epsilon \dot{u}_0(s)}_{(H^1)^*}^2 \di s & \leq \norm{\Rop_\epsilon}_{\textnormal{op}}^2 \int_0^t \norm{\dot{u}_{0}(s)}^2_{H^2}\di s \leq C\epsilon \int_0^t \norm{\dot{u}_{0}(s)}^2_{H^2}\di s,
\end{aligned}
\]
and \cref{eq:temp per 3} boils down to
\[
\norm{w_\epsilon(t)}_{H^1}^2 \leq C\Big\{\epsilon\norm{u_{0}}^2_{L^\infty(I;H^2)}+\epsilon \int_0^t \norm{\dot{u}_{0}(s)}^2_{H^2}\di s+ \int_0^t\norm{w_\epsilon(s)}_{H^1}^2\di s  \Big\},
\]
where the positive constant $C$ only depends on $m$, $M$ and $\Omega$. Finally, we complete the proof by Gr\"{o}nwall's inequality and the embedding $H^1(I;H^2)\hookrightarrow L^\infty(I;H^2)$.
\end{proof}

Now we turn to prove \cref{thm:main esti}.
\begin{proof}
Let $w_\epsilon$ be defined as \cref{lem:mid esti}. In order to estimate $w_\epsilon(t)-r_\epsilon(t)$, we are left to handle
\[
\epsilon \chi_j^{\alpha\beta}(x/\epsilon)\theta_\epsilon \Sop_\epsilon^2\RoundBrackets*{\partial_j u_0^\beta(t)}-\epsilon \chi_j^{\alpha\beta}(x/\epsilon)\partial_j u_0^\beta(t).
\]
Note that
\[
\begin{aligned}
&\partial_k \RoundBrackets*{\epsilon \chi_j^{\alpha\beta}\RoundBrackets*{x/\epsilon}\theta_\epsilon \Sop_\epsilon^2\RoundBrackets*{\partial_j u_0^\beta(t)}}\\
=& \RoundBrackets*{\partial_k \chi_j^{\alpha\beta}}\RoundBrackets*{x/\epsilon}\theta_\epsilon \Sop_\epsilon^2\RoundBrackets*{\partial_j u_0^\beta(t)} \\
&+\chi_j^{\alpha\beta}\RoundBrackets*{x/\epsilon} \partial_k \RoundBrackets*{\epsilon\theta_\epsilon} \Sop_\epsilon^2\RoundBrackets*{\partial_j u_0^\beta(t)} \\
&+\epsilon \chi_j^{\alpha\beta}\RoundBrackets*{x/\epsilon}\theta_\epsilon \Sop_\epsilon^2\RoundBrackets*{\partial_{jk} u_0^\beta(t)} \\
\coloneqq& J_1+J_2+J_3,
\end{aligned}
\]
where the property $\partial_j \Sop_\epsilon(f)=\Sop_\epsilon(\partial_jf)$ is utilized (ref. \cite{Shen2017,Shen2018}). Similarly,
\[
\begin{aligned}
\partial_k\RoundBrackets*{\epsilon \chi_j^{\alpha\beta}\RoundBrackets*{x/\epsilon} \partial_j u_0^\beta(t)}=& \RoundBrackets*{\partial_k \chi_j^{\alpha\beta}}\RoundBrackets*{x/\epsilon}\partial_j u_0^\beta(t)+\epsilon\chi_j^{\alpha\beta}\RoundBrackets*{x/\epsilon}\partial_{jk} u_0^\beta(t)\\
\coloneqq& K_1+K_2.
\end{aligned}
\]
According to the boundary layer estimate (i.e. Proposition 3.1.8. of \cite{Shen2018}), we will have
\[
\norm{J_2}_{L^2(\Omega)} \leq C\epsilon^{1/2} \max_{j,\beta}\norm{\chi_j^{\alpha\beta}}_{L^2(Y)} \norm{u_{0}(t)}_{H^2}.
\]
For $J_3$, we can show that
\[
\norm{J_3}_{L^2(\Omega)} \leq C\epsilon \max_{j,\beta}\norm{\chi_j^{\alpha\beta}}_{L^2(Y)} \norm{u_{0}(t)}_{H^2}
\]
by Proposition 3.1.5. of \cite{Shen2018}. The estimates for $K_2$ will be straightforward:
\[
\norm{K_2}_{L^2(\Omega)} \leq \epsilon \max_{j,\beta}\norm{\chi_j^{\alpha\beta}}_{L^\infty(Y)} \norm{u_{0}(t)}_{H^2}.
\]
For $J_1-K_1$, we can show that
\[
\begin{aligned}
&\norm{J_1-K_1}_{L^2(\Omega)} \\
\leq & \norm{ \RoundBrackets*{\partial_k \chi_j^{\alpha\beta}}\RoundBrackets*{x/\epsilon}\theta_\epsilon \Sop_\epsilon^2\RoundBrackets*{\partial_j u_0^\beta(t)} -  \RoundBrackets*{\partial_k\chi_j^{\alpha\beta}}\RoundBrackets*{x/\epsilon} \Sop_\epsilon^2\RoundBrackets*{\partial_j u_0^\beta(t)}}_{L^2(\Omega)} \\
&+ \norm{ \RoundBrackets*{\partial_k\chi_j^{\alpha\beta}}\RoundBrackets*{x/\epsilon} \Sop_\epsilon^2\RoundBrackets*{\partial_j u_0^\beta(t)} - \RoundBrackets*{\partial_k\chi_j^{\alpha\beta}}\RoundBrackets*{x/\epsilon} \partial_j u_0^\beta(t)}_{L^2(\Omega)} \\
\coloneqq& N_1+N_2,
\end{aligned}
\]
where $N_1 \leq C\epsilon^{1/2} \sum_{j,\beta}\norm{\partial_k\chi_j^{\alpha\beta}}_{L^2(\Omega)} \norm{u_{0}(t)}_{H^2}$ is again from Proposition 3.1.8. of \cite{Shen2018} and
\[
\begin{aligned}
N_2 \leq & \sum_{j,\beta} \norm{\partial_k \chi_j^{\alpha\beta}}_{L^\infty(Y)}\norm{ \Sop_\epsilon^2\RoundBrackets*{\partial_j u_0^\beta(t)} - \partial_j u_0^\beta(t) }_{L^2(\Omega)} \\
\leq & C \epsilon \max_{k,\gamma} \norm{\chi_k^\gamma}_{W^{1,\infty}_\#(Y;\Real^d)}  \norm{u_{0}(t)}_{H^2},
\end{aligned}
\]
here in the last line above we have used Proposition 3.1.6. of \cite{Shen2018}. Furthermore, the estimate of
\[
\norm{w_\epsilon(t)-r_\epsilon(t)}_{L^2(\Omega;\Real^d)} \leq C\epsilon\CurlyBrackets*{\sum_{j,\beta} \RoundBrackets*{\norm{\chi_j^{\alpha\beta}}_{L^2(\Omega)}+\norm{\chi_j^{\alpha\beta}}_{L^\infty(\Omega)}}\norm{\partial_ju_{0}^\beta(t)}_{L^2(\Omega)}}
\]
is straightforward. Thus, combining all the estimates above, we derive that
\[
\norm{w_\epsilon(t)-r_\epsilon(t)}_{H^1} \leq C\epsilon^{1/2} \max_{k,\gamma} \norm{\chi_k^\gamma}_{W^{1,\infty}_\#(Y;\Real^d)} \norm{u_{0}(t)}_{H^2},
\]
which completes the proof.
\end{proof}

\Section{Numerical experiments} \label{sec:Num}

The numerical experiments in this section are carried out to validate the convergence estimate obtained in \cref{thm:main esti}. We consider a 2D square domain $\Omega=(0,1)\times(0,1)$ and equally divide it into $N\times N$ entirely periodic cells, means that $\epsilon=1/N$. Each cell is constituted by two patterns of isotropic linear elastic materials, i.e., reinforced phase and matrix, we denote the pairs of Young's modulus and Possion's ratio by $(E_0,\nu_0)$ and $(E_1,\nu_1)$, respectively. In \emph{plane stress problems} (see \cite{Sadd2021}), the elastic tensor $\tenA$ could be expressed by $(E,\nu)$ as $\tenA_{ij}^{\alpha\beta} = \lambda \delta_{\alpha i} \delta_{\beta j} + \mu\RoundBrackets{\delta_{\alpha \beta}\delta_{ij}+\delta_{\alpha j}\delta_{\beta i}}$, where
\[
\lambda = \frac{E\nu}{(1+\nu)(1-\nu)} \  \mathup{and} \  \mu = \frac{E}{2(1+\nu)},
\]
noticing here $\lambda$ differs from the first Lam\'{e} parameter in 3D elastic problems.

\begin{figure}[htbp]
\centering
\begin{subfigure}[b]{0.49\textwidth}
\centering
\begin{tikzpicture}
\def\xa{0};
\def\xb{2};
\def\ya{0};
\def\yb{2};
\coordinate (Oy) at ({\xa-1}, {\ya-1});
\coordinate (Ry0) at (\xa, \ya);
\coordinate (Ry1) at (\xb, \ya);
\coordinate (Ry2) at (\xa, \yb);
\coordinate (Ry3) at (\xb, \yb);

\draw[axis] (Oy) -- ++(2, 0) node[below] {$y_1$};
\draw[axis] (Oy) -- ++(0, 2) node[left] {$y_2$};

\fill[GeoDataViz-D3] (Ry0) rectangle (Ry3);

\fill[GeoDataViz-D4] ({7/8*\xa+1/8*\xb}, {5/8*\ya+3/8*\yb}) rectangle ({1/8*\xa+7/8*\xb}, {3/8*\ya+5/8*\yb});
\fill[GeoDataViz-D4] ({5/8*\xa+3/8*\xb}, {7/8*\ya+1/8*\yb}) rectangle ({3/8*\xa+5/8*\xb}, {1/8*\ya+7/8*\yb});
\draw[label] ({0.5*\xa+0.5*\xb}, {0.5*\ya+0.5*\yb}) -- ++(-2.5, 1.5) node[above] {\small $(E_0, \nu_0)$};
\draw[label] ({0.75*\xa+0.25*\xb}, {0.75*\ya+0.25*\yb}) -- ++(-0.75, -0.75) node[below] {\small ${(E_1, \nu_1)}$};
\draw[dotted] ({\xa-0.25}, {5/8*\ya+3/8*\yb}) -- ({\xb+0.25}, {5/8*\ya+3/8*\yb}) node[right] {\tiny $3/8$};
\draw[dotted] ({\xa-0.25}, {7/8*\ya+1/8*\yb}) -- ({\xb+0.25}, {7/8*\ya+1/8*\yb}) node[right] {\tiny $1/8$};
\draw[dotted] ({\xa-0.25}, {3/8*\ya+5/8*\yb}) -- ({\xb+0.25}, {3/8*\ya+5/8*\yb}) node[right] {\tiny $5/8$};
\draw[dotted] ({\xa-0.25}, {1/8*\ya+7/8*\yb}) -- ({\xb+0.25}, {1/8*\ya+7/8*\yb}) node[right] {\tiny $7/8$};

\draw[dotted] ({5/8*\xa+3/8*\xb}, {\ya-0.25}) -- ({5/8*\xa+3/8*\xb}, {\yb+0.25}) node[above] {\tiny $\frac{3}{8}$};
\draw[dotted] ({3/8*\xa+5/8*\xb}, {\ya-0.25}) -- ({3/8*\xa+5/8*\xb}, {\yb+0.25}) node[above] {\tiny $\frac{5}{8}$};
\draw[dotted] ({7/8*\xa+1/8*\xb}, {\ya-0.25}) -- ({7/8*\xa+1/8*\xb}, {\yb+0.25}) node[above] {\tiny $\frac{1}{8}$};
\draw[dotted] ({1/8*\xa+7/8*\xb}, {\ya-0.25}) -- ({1/8*\xa+7/8*\xb}, {\yb+0.25}) node[above] {\tiny $\frac{7}{8}$};

\end{tikzpicture}
\caption{}\label{fig:cell-a}
\end{subfigure}
\begin{subfigure}[b]{0.49\textwidth}
\centering
\begin{tikzpicture}
\def\xa{0};
\def\xb{2};
\def\ya{0};
\def\yb{2};
\coordinate (Oy) at ({\xa-1}, {\ya-1});
\coordinate (Ry0) at (\xa, \ya);
\coordinate (Ry1) at (\xb, \ya);
\coordinate (Ry2) at (\xa, \yb);
\coordinate (Ry3) at (\xb, \yb);

\draw[axis] (Oy) -- ++(2, 0) node[below] {$y_1$};
\draw[axis] (Oy) -- ++(0, 2) node[left] {$y_2$};

\fill[GeoDataViz-D3] (Ry0) rectangle (Ry3);

\fill[GeoDataViz-D4] ({\xa}, {5/8*\ya+3/8*\yb}) rectangle ({\xb}, {3/8*\ya+5/8*\yb});
\fill[GeoDataViz-D4] ({5/8*\xa+3/8*\xb}, {\ya}) rectangle ({3/8*\xa+5/8*\xb}, {\yb});
\draw[label] ({0.5*\xa+0.5*\xb}, {0.5*\ya+0.5*\yb}) -- ++(-2.5, 1.5) node[above] {\small $(E_0, \nu_0)$};
\draw[label] ({0.75*\xa+0.25*\xb}, {0.75*\ya+0.25*\yb}) -- ++(-0.75, -0.75) node[below] {\small ${(E_1, \nu_1)}$};
\draw[dotted] ({\xa-0.25}, {5/8*\ya+3/8*\yb}) -- ({\xb+0.25}, {5/8*\ya+3/8*\yb}) node[right] {\tiny $3/8$};
\draw[dotted] ({\xa-0.25}, {3/8*\ya+5/8*\yb}) -- ({\xb+0.25}, {3/8*\ya+5/8*\yb}) node[right] {\tiny $5/8$};
\draw[dotted] ({5/8*\xa+3/8*\xb}, {\ya-0.25}) -- ({5/8*\xa+3/8*\xb}, {\yb+0.25}) node[above] {\tiny $\frac{3}{8}$};
\draw[dotted] ({3/8*\xa+5/8*\xb}, {\ya-0.25}) -- ({3/8*\xa+5/8*\xb}, {\yb+0.25}) node[above] {\tiny $\frac{5}{8}$};
\end{tikzpicture}
\caption{}\label{fig:cell-b}
\end{subfigure}
\begin{subfigure}[b]{0.49\textwidth}
\centering
\begin{tikzpicture}[scale=0.95]
\def\xa{0};
\def\xb{4};
\def\ya{0};
\def\yb{4};

\coordinate (Ox) at (\xa, \ya);
\coordinate (A) at (\xb, \ya);
\coordinate (B) at (\xa, \yb);
\coordinate (C) at (\xb, \yb);

\draw[axis] (Ox) -- ++({\xb+1}, \ya) node[below] {$x_1$};
\draw[axis] (Ox) -- ++(\xa, {\yb+1}) node[left] {$x_2$};

\fill[GeoDataViz-D3] (Ox) rectangle (C);
\foreach \x in {0.00, 0.25, 0.5, 0.75}
\foreach \y in {0.00, 0.25, 0.5, 0.75}
{
\def\xaCell{\xb*\x+\xa*(1.0-\x)};
\def\xbCell{\xb*(\x+0.25)+\xa*(0.75-\x)};
\def\yaCell{\yb*\y+\ya*(1.0-\y)};
\def\ybCell{\yb*(\y+0.25)+\ya*(0.75-\y)};
\fill[GeoDataViz-D4] ({7/8*\xaCell+1/8*\xbCell}, {5/8*\yaCell+3/8*\ybCell}) rectangle ({1/8*\xaCell+7/8*\xbCell}, {3/8*\yaCell+5/8*\ybCell});
\fill[GeoDataViz-D4] ({5/8*\xaCell+3/8*\xbCell}, {7/8*\yaCell+1/8*\ybCell}) rectangle ({3/8*\xaCell+5/8*\xbCell}, {1/8*\yaCell+7/8*\ybCell});
%\fill[GeoDataViz-D4] ({\xaCell}, {5/8*\yaCell+3/8*\ybCell}) rectangle ({\xbCell}, {3/8*\yaCell+5/8*\ybCell});
%\fill[GeoDataViz-D4] ({5/8*\xaCell+3/8*\xbCell}, {\yaCell}) rectangle ({3/8*\xaCell+5/8*\xbCell}, {\ybCell});
}

\def\TriH{0.4};
\foreach \x in {1,...,6}
{
\def\Triy{\x/7*\yb};
\filldraw[fill=gray] (\xb, \Triy) -- ++(\TriH, {0.5*\TriH}) -- ++(0, {-\TriH});
}

\filldraw[pattern=north west lines,pattern color=GeoDataViz-D3] ({\xa-0.5}, \ya) rectangle ({\xb+0.5}, {\ya-0.5});
\draw[load] ({0.5*\xa+0.5*\xb}, {0.5*\ya+0.5*\yb}) -- ++(0, -1.5) node[right] {\small $f$};

\foreach \x in {1,...,6}
{
\def\Tracy{\x/7*\yb};
\draw[load] ({\x/7-1.25}, {\x/7*\yb+0.125}) -- (\xa, \Tracy);
}

\node[text width=1em] at (-1, 2) {\small $t$};
\end{tikzpicture}
\caption{}\label{fig:demo BCs}
\end{subfigure}

\caption{(a) and (b), two cell configurations; (c) an illustration of the domain with the cell configuration (a), where $N=4$ and a homogeneous Dirichlet boundary condition is set on $\{1\}\times(0,1)$, a boundary traction is applied on $\{0\}\times(0,1)$, Tresca's law is modeled on $(0,1)\times\{0\}$.}
\end{figure}
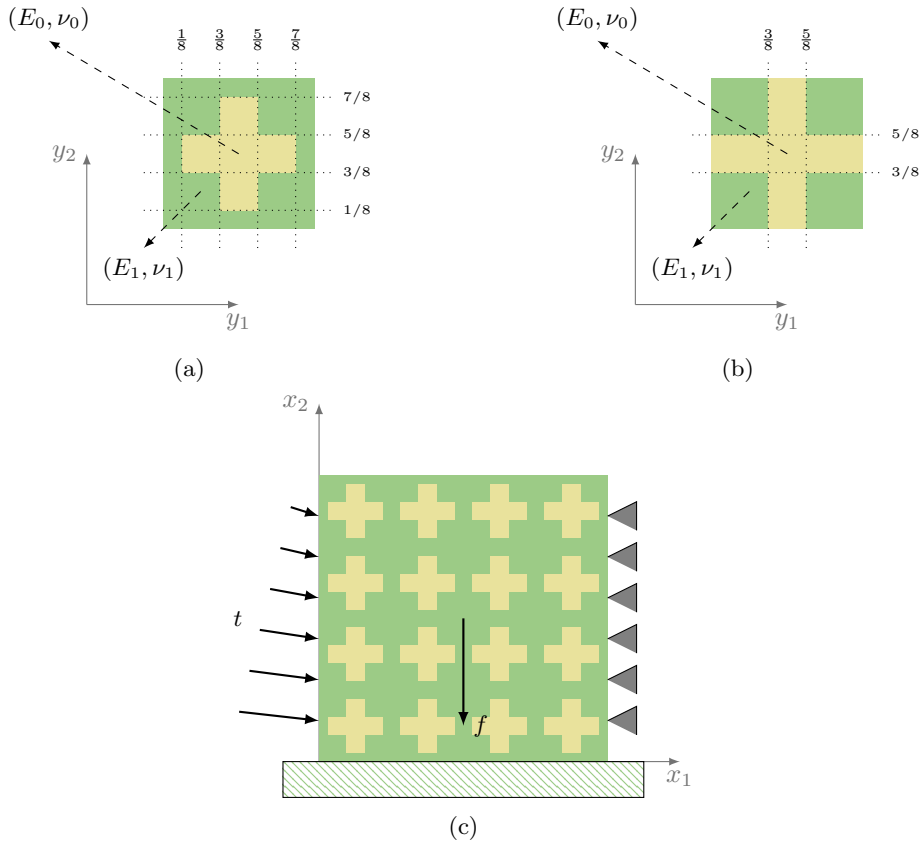

Our experiments are for the cell configurations in \cref{fig:cell-a,fig:cell-b}, which determine $1$-periodic tensor-valued functions $\tenA(y)$. We consider the problem \cref{eq:Tresca weak}, while the elastic tensor $\tenA$ is replaced by $\tenA(x/\epsilon)$ (the original small periodic problem) or $\hat{\tenA}$ (the homogenized problem), and the settings of boundary conditions are demonstrated in \cref{fig:demo BCs}. All parameter values and solvability conditions are listed in \cref{tab:parameters}. Note that \cref{fig:demo BCs} demonstrates the elastic body that is formed by periodically duplicating the cell \cref{fig:cell-a}, and only the matrix corresponding to the material parameters $(E_1, \nu_1)$ contacts with the rigid body (see the lower boundary of \cref{fig:demo BCs} for the illustration). Therefore, the value of $\TrescaBound$ in \cref{tab:parameters} has a solid physical background and can be measured from a physical experiment. However, if the domain consists of the cells of \cref{fig:cell-b}, then $\TrescaBound$ is ambiguous because reinforced phase also contacts on $\GammaC$. In certain situations, according the force balance law on the whole domain, we may obtain the total internal force applied on the contact boundary by conducting physical experiments. By using the averaged friction coefficient, the value of $\TrescaBound$ could also be estimated by the product of the averaged friction coefficient and normal internal force.

\begin{table}[ht]
\caption{The parameter values and solvability conditions used in numerical experiments.}\label{tab:parameters}
\centering
{\begin{tabular}{@{}cc@{}}
\hline
$(E_0, \nu_0), (E_1, \nu_1)$ & $(\qty{77.2}{GPa}, 0.33), (\qty{117.0}{GPa}, 0.43)$ \\
%$(E_0, \nu_0), (E_1, \nu_1)$ & $({77.2}\mathup{ GPa}, 0.33), ({117.0}\mathup{ GPa}, 0.43)$ \\ \hline
\hline
$T$ & \qty{1.0}{s} \\
%$T$ & ${1.0}\mathup{ s}$ \\ \hline
\hline
$f(x_1,x_2,t)$ & $(0.0, \num{-1.0e-4})\,\mathup{GN/m}^3$ \\
\hline
$ t(x_1,x_2,t)$ & $(0.08(1.25-x_2)t, -0.01t)\,\mathup{GPa}$ \\
\hline
$\TrescaBound$ & $\qty{0.004}{GPa}$ \\
\hline
%$H_\textnormal{T}$ & ${0.004}\mathup{ GPa}$ \\ \hline
\end{tabular}}
\end{table}

According to the assumptions \cref{eq:init k per,eq:init inf per}, we have $u_{\epsilon,\star}(x)=u_{0,\star}(x) \equiv 0$ for the initial conditions. To solve \cref{eq:Tresca weak}, we first take a semi-discretization of $u$ w.r.t. the temporal variable, that is $\dot{u}(t) \leftarrow \delta u^m \coloneqq \RoundBrackets{u^m-u^{m-1}}/\Delta t$, where $\Delta t$ is the time step size. We then convert \cref{eq:Tresca weak} into a variational inequality of the second kind (see \cite{Han1999}) for $\delta u^m$ , which is corresponding to the minimizer of a nonsmooth convex functional. Some previous efforts from the author on solving contact problems could be found in \cite{Cui1983,Cui1988}, and interested readers may refer to \cite{Cui1988} for a comparative study of numerical methods with real observed data.

In our numerical experiments, we partition every cell into $32\times32$ elements, which says that the original problem and homogenized problem are solved on a $1024\times 1024$ mesh for $N=32$, and we discretize the whole system in the bilinear finite element space. As for temporal discretizations, we fix the time step size as $\Delta t= 1/1024$ for all numerical cases. We implement a Nestrov accelerated proximal gradient descent method (ref. \cite{Nesterov2012}) to solve the nonsmooth optimization problem of every time step. Two notations are introduced to measure asymptotic errors:
\[
\begin{aligned}
%\textbf{Err}_0 &\coloneqq \sup_{t\in (0,T]} \frac{\norm{u_\epsilon(t)-u_{0}(t)}_{L^2}}{\norm{u_{0}(t)}_{L^2}}, \\
\textbf{Err}_1 &\coloneqq \sup_{t\in (0,T]} \frac{\norm{\nabla \big(u_\epsilon(t)-u_{0}(t)\big)}_{L^2}}{\norm{\nabla u_{0}(t)}_{L^2}}, \\
\textbf{Err}_2 &\coloneqq \sup_{t\in (0,T]} \frac{\CurlyBrackets*{\sum_{i,\alpha}\norm{\partial_i u_\epsilon^\alpha(t) -\partial_i \RoundBrackets*{u_0^\alpha(t) +\epsilon \chi_j^{\alpha \beta}(x/\epsilon) \partial_j u_{0}^\beta(t) }  }^2_{L^2}}^{1/2}}{\norm{\nabla u_{0}(t)}_{L^2}},
\end{aligned}
\]
where $u_{\epsilon}$ and $u_{0}$ are numerical solutions. Since $u_{0}$ obtained by numerical method are only piecewisely bilinear, taking partial derivatives to order $2$ is problematic in $\textbf{Err}_2$. Therefore, calculations of $L^2$-norm in $\textbf{Err}_2$ are first conducted on every element then taking a summation. The numerical results are listed in Table 2, and our algorithms are build up from PETSc (\cite{Balay2022}) and source codes are hosted on GitHub\footnote{https://github.com/Laphet/MS-TrescaBP.git}.

\begin{table}[htbp]
\centering
{\footnotesize
\caption{(a) Numerical results for the cell configuration in \cref{fig:cell-a}, (b) for the cell configuration in \cref{fig:cell-b}.} \label{tab:results-abcdefgh}
\begin{subtable}[b]{0.49\textwidth}
\caption{}\label{tab:data-abcd}
\centering
\begin{tabular}{@{}ccccc@{}}
\hline
 $N$ & $\epsilon$ & $h $ & $\textbf{Err}_1$ & $\textbf{Err}_2$ \\
\hline
$4$ & $1/4$ & $1/128$ & $0.17981$ & $0.05691$ \\
$8$ & $1/8$ & $1/256$ & $0.18118$ & $0.04087$ \\
$16$ & $1/16$ & $1/512$ & $0.18116$ & $0.02897$ \\
$32$ & $1/32$ & $1/1024$ & $0.18186$ & $0.02043$ \\
\hline
\end{tabular}
\end{subtable}
\begin{subtable}[b]{0.49\textwidth}
\caption{}\label{tab:data-efgh}
\centering
\begin{tabular}{@{}ccccc@{}}
\hline
 $N$ & $\epsilon$ & $h$ & $\textbf{Err}_1$ & $\textbf{Err}_2$ \\
\hline
$4$ & $1/4$ & $1/128$ & $0.20064$ & $0.08512$ \\
$8$ & $1/8$ & $1/256$ & $0.19946$ & $0.06018$ \\
$16$ & $1/16$ & $1/512$ & $0.19910$ & $0.04240$ \\
$32$ & $1/32$ & $1/1024$ & $0.19914$ & $0.02987$ \\
\hline
\end{tabular}
\end{subtable}
}
\end{table}

We can see that $\textbf{Err}_2 \approx O(\sqrt{\epsilon})$, which agrees satisfactorily with our theorem. Meanwhile, those data of $\textbf{Err}_1$ show that $u_{0}$ cannot approximate $u_\epsilon$ in gradients. %The values of $\textbf{Err}_0$ in \cref{tab:data-efgh} imply that $\norm{u_{\epsilon}-u_{0}}_{L^2(\Omega)} \approx O(\epsilon)$, while a record ($N=16$) in \cref{tab:data-abcd} tends to negative this conjecture.
In general, it will be much harder to prove or disprove an optimal estimate of asymptotic errors in $L^2$-norm, and existing methods are mostly based on the maximum principle or the duality technique, which are both inapplicable for our problems.

\Section{Conclusions and further discussions}\label{sec:Conclu}

The study of this paper originates from the modeling of friction contacts of composite material structures. Specifically, we consider the quasistatic Tresca friction law, which could be categorized as non-smooth and nonlinear contact BC. Our theoretical results are based on two different homogenization settings. One is H-convergences, that the coefficients of the linear elasticity system H-converge to the homogenized ones, and we obtain a qualitative theorem which states that the oscillating solutions converge weakly to the homogenized solution in a proper space. Another is the homogenization with small periodicity, it means that the coefficients have small periodicity, and then we derive a quantitative theorem on the estimate of asymptotic errors in $H^1$-norm. Our numerical experiments show that the convergence rate seems to be optimal.

In our numerical illustrations, the Tresca friction bound is independent of the heterogeneity of composite materials, which is an approximation to the friction experiment setting for the structural systems of composite materials. Here we assume that the reinforced phases (e.g., $(E_0,\nu_0)$ in \cref{fig:cell-a}) have a positive distance away from contact boundaries and the friction contacts only happen on matrix materials (e.g., $(E_1,\nu_1)$ in \cref{fig:cell-a}). It should be admitted that allowing material-dependent BCs is challenging in mathematical analysis, and existing results show that this problem is related to some delicate descriptions of the domain (ref. \cite{GerardVaret2012,Shen2020,Geng2020}).

It is worth pointing out, in real problems simulated, the homogenized solution may not meet the regularity assumption of \cref{thm:main esti}. In this case, the first-order expression may be a poor approximation to the original oscillating solution, and could not precisely capture the refined information on the contact boundaries. However, the Newton-like methods are commonly adopted to solve frictional contact problems, while its performance is hinged on linear solvers for each Newton step, which could be accelerated by using first-order expressions as initial guesses.

The Coulomb friction contact problems of the structural systems made from composite materials essentially have the multiscale and strongly nonlinear features, especially on the contact boundaries, and it is difficult to analyze them in mathematics. Another alternative model for Coulomb’s law is the so-called normal compliance model (\cite{Klarbring1988,Shillor2004}), which in some sense is friendly to mathematical analysis. Hence, it is natural to consider the homogenization with normal compliance contact laws which will be our future work.

\bigskip

%{\bf Acknowledgement}\ \ The content of the Acknowledgement.

\bibliographystyle{plain}
\bibliography{refs}
%\begin{thebibliography}{aa}
%
%\bibitem{1}%[1]
%Jacobson, N. and Smith, J., Lie  Algebras, Dover., Publ., New York,
%1979.
%
%\bibitem{L}%[2]
%Chew, B. S., On the commutativity of restricted Lie algebras, {\it
%Proc. Amer. Math. Soc}, {\bf 16}(3), 1965, 547--560.
%
%\bibitem{C}%[10]
%Feldvoss, J., Projective modules and block  of supersolvable
%restricted Lie algebras, {\it J. Algebra,} {\bf 222}, 1999,
%284--300.
%
%\end{thebibliography}

\end{document}